\documentclass[12pt]{amsart}
\input{amssymb.sty}
\title{Solid von Neumann Algebras}
\author{Narutaka OZAWA}
\address{Department of Mathematical Sciences,
University of Tokyo, Komaba, 153-8914}
\email{narutaka@ms.u-tokyo.ac.jp}
\date{Feburuary 7, 2003}
\thanks{The author was supported by the JSPS 
Postdoctoral Fellowships for Research Abroad}

\dedicatory{Dedicated to Professor Masamichi Takesaki 
on the occasion of his 70th birthday.}
\newtheorem{thm}{Theorem}
\newtheorem{prop}[thm]{Proposition}
\newtheorem{lem}[thm]{Lemma}
\newcommand{\G}{\Gamma} 
\newcommand{\e}{\varepsilon} 
\newcommand{\B}{{\mathbb B}}  
\newcommand{\K}{{\mathbb K}}  
\newcommand{\F}{{\mathbb F}} 
\newcommand{\EE}{{\mathcal E}} 
\newcommand{\FF}{{\mathcal F}} 
\newcommand{\LG}{{\mathcal L}\G} 
\newcommand{\MM}{{\mathcal M}} 
\newcommand{\NN}{{\mathcal N}} 
\newcommand{\A}{{\mathcal A}} 
\newcommand{\BB}{{\mathcal B}} 
\newcommand{\ce}{\e_{\NN}} 
\newcommand{\hh}{{\mathcal H}} 
\newcommand{\eg}{\ell_2\G} 
\newcommand{\id}{\mathrm{id}} 
\newcommand{\clg}{C^*_\lambda\G}
\newcommand{\crg}{C^*_\rho\G}
\newcommand{\ip}[1]{\langle#1\rangle} 
\begin{document}
\begin{abstract}
We prove that the relative commutant of 
a diffuse von Neumann subalgebra in 
a hyperbolic group von Neumann algebra 
is always injective. 
It follows that any non-injective subfactor 
in a hyperbolic group von Neumann algebra 
is non-$\Gamma$ and prime. 
The proof is based on $C^*$-algebra theory. 
\end{abstract}
\maketitle
\section{Introduction}
Recall that a von Neumann algebra is said to be diffuse 
if it does not contain a minimal projection. 
We say a von Neumann algebra $\MM$ is \textit{solid} 
if for any diffuse von Neumann subalgebra $\A$ in $\MM$, 
the relative commutant $\A'\cap\MM$ is injective. 
A solid von Neumann algebra is necessarily finite. 
We prove the following theorem which answers a question of Ge \cite{ge} 
that whether the free group factors are solid. 

\begin{thm}\label{thm} 
The group von Neumann algebra $\LG$ of 
a hyperbolic group $\G$ is solid. 
\end{thm}

Recall that a factor $\MM$ is said to be prime 
if $\MM\cong\MM_1\bar{\otimes}\MM_2$ implies 
either $\MM_1$ or $\MM_2$ is finite dimensional. 
The existence of such factors was proved by Popa \cite{popa} 
who showed the group factors of uncountable free groups are prime. 
The case for countable free groups had remained open for some time, 
but was settled by Ge \cite{ge}. 
This was generalized by \c{S}tefan \cite{stefan} to 
their subfactors of finite index.
Our theorem gives a further generalization. 
Indeed, combined with a result of Popa \cite{popapc} 
(Proposition \ref{popa} in this paper), 
we obtain the following proposition. 

\begin{prop}\label{prop}
A subfactor of a solid factor is again solid and 
a solid factor is non-$\G$ and prime unless it is injective. 
\end{prop}

Notably, this provides infinitely many prime $\mathrm{II}_1$-factors 
(with the property $\mathrm{(T)}$). 
Indeed thanks to a theorem of Cowling and Haagerup \cite{ch}, 
for lattices $\G_n$ in $\mathrm{Sp}(1,n)$, 
we have $\LG_m\not\cong \LG_n$ whenever $m\neq n$. 
However, we are unaware of any non-injective solid factor 
with the Haagerup property other than the free group factor(s). 
This proposition also distinguishes the factor 
$({\mathcal L}\F_\infty\bar{\otimes}L_\infty[0,1])
\ast{\mathcal L}\F_\infty$ 
from the free group factor ${\mathcal L}\F_\infty$, 
which answers a question of Shlyakhtenko. 
\section{Preliminary Results on Reduced Group $C^*$-algebras}\label{grp}

For a discrete group $\G$,
we denote by $\lambda$ (resp.\ $\rho$) the left (resp.\ right) 
regular representation on $\eg$ and let 
$\clg$ (resp. $\crg$) be the $C^*$-algebra generated 
by $\lambda(\G)$ (resp.\ $\rho(\G)$) in $\B(\eg)$ and 
$\LG=(\clg)''$ be its weak closure. 
The $C^*$-algebra $\clg$ (resp. the von Neumann algebra $\LG$) 
is called the reduced group $C^*$-algebra 
(resp.\ the group von Neumann algebra). 

The study of $C^*$-norms on tensor products was initiated 
in the 50's by Turumaru, but the first substantial result 
was obtained by Takesaki \cite{takesaki} who showed 
the minimal tensor norm is the smallest among 
the possible $C^*$-norms on a tensor product of $C^*$-algebras. 
He also introduced the notion of nuclearity and 
found that the reduced group $C^*$-algebra $C^*_\lambda\F_2$ of
the free group $\F_2$ on two generators is not nuclear. 
Namely, the $*$-homomorphism 
\[
C^*_\lambda\F_2\otimes C^*_\rho\F_2\ni\sum_{i=1}^na_i\otimes x_i
\mapsto \sum_{i=1}^na_ix_i\in\B(\ell_2\F_2)
\] 
is not continuous w.r.t.\ the minimal tensor norm. 
Yet, Akemann and Ostrand \cite{ao} proved a remarkable theorem
that it is continuous if one composes it with the quotient map 
$\pi$ from $\B(\ell_2\F_2)$ onto 
the Calkin algebra $\B(\ell_2\F_2)/\K(\ell_2\F_2)$. 
By a completely different argument, 
Skandalis (Th\'{e}or\`{e}me 4.4 in \cite{skandalis}) 
proved the same for all discrete subgroups 
in connected simple Lie groups of rank one, 
and Higson and Guentner (Lemma 6.2.8 in \cite{hg}) 
for all hyperbolic groups.
In summary, 

\begin{thm}\label{ao}
Let $\G$ be a hyperbolic group 
or a discrete subgroup in a connected simple Lie group of rank one. 
Then, the $*$-homomorphism 
\[
\nu_{\G}\colon\clg\otimes\crg\ni\sum_{i=1}^na_i\otimes x_i
\mapsto\pi(\sum_{i=1}^na_ix_i)\in\B(\eg)/\K(\eg)
\]
is continuous w.r.t.\ the minimal tensor norm on $\clg\otimes\crg$. 
\end{thm}

The crucial ingredient in the proof was the amenability 
of the action of $\G$ on 
a suitable boundary which is `small at infinity'. 
For the information on amenable actions, 
we refer the reader to the book \cite{ar} of 
Anantharaman-Delaroche and Renault.
Since $\G$ acts amenably on a compact set, 
$\clg$ is embeddable into a nuclear $C^*$-algebra and thus has 
the property (C) of Archbold and Batty (Theorem 3.6 in \cite{ab}). 
Although we do not need this fact, we mention that 
the property (C) is equivalent to exactness 
by a deep theorem of Kirchberg \cite{kirchberg}.
By Effros and Haagerup's theorem (Theorem 5.1 in \cite{eh}), 
the property (C) implies the local reflexivity. 
In summary, 
\begin{lem}\label{locref}
Let $\G$ be as above. 
Then, $\clg$ is locally reflexive, i.e., 
for any finite dimensional operator system $E\subset(\clg)^{**}$, 
there is a net of unital completely positive maps 
$\theta_i\colon E\to\clg$ which converges to $\id_E$ 
in the point-weak$^*$ topology. 
\end{lem}

\section{Proof of The Theorem}
We recall the following principle \cite{choi}; 
if $\Psi\colon A\to B$ is a unital completely positive map 
and its restriction to a $C^*$-subalgebra $A_0\subset A$ 
is multiplicative, 
then we have $\Psi(axb)=\Psi(a)\Psi(x)\Psi(b)$ 
for any $a,b\in A_0$ and $x\in A$. 

Let $\NN\subset\MM$ be finite von Neumann algebras 
with a faithful trace $\tau$ on $\MM$. 
Then, there is a normal conditional expectation 
$\ce$ from $\MM$ onto $\NN$, which is defined by the relation 
$\tau(\ce(a)x)=\tau(ax)$ for $a\in\MM$ and $x\in\NN$. 
This implies that a von Neumann subalgebra of 
a finite injective von Neumann algebra is again injective. 
Moreover, $\ce$ is unique in the sense that 
any trace preserving conditional expectation 
from $\MM$ onto $\NN$ coincides with $\ce$. 
Indeed, for any $a\in\MM$ and $x\in\NN$, we have
\[
\tau(\e'(a)x)=\tau(\e'(ax))=\tau(ax)=\tau(\ce(a)x).
\]

We say a von Neumann subalgebra $\MM$ in $\B(\hh)$ 
satisfies the condition (AO) if 
\begin{quote}
there are unital ultraweakly dense $C^*$-subalgebras 
$B\subset\MM$ and $C\subset\MM'$ such that 
$B$ is locally reflexive and the $*$-homomorphism 
\[
\nu\colon B\otimes C\ni\sum_{i=1}^na_i\otimes x_i
\mapsto\pi(\sum_{i=1}^na_ix_i)\in\B(\hh)/\K(\hh)
\]
is continuous w.r.t.\ the minimal tensor norm on $B\otimes C$. 
\end{quote}
We have seen in Section \ref{grp} that the group von Neumann algebra $\LG$ 
satisfies the condition (AO) whenever $\G$ is a hyperbolic group 
or a discrete subgroup in a connected simple Lie group of rank one. 

\begin{lem}\label{mc}
Let $B\subset\MM$ and $C\subset\MM'$ 
be unital ultraweakly dense $C^*$-subalgebras 
with $B$ locally reflexive and 
let $\NN\subset\MM$ be a von Neumann subalgebra 
with a normal conditional expectation $\ce$ onto $\NN$. 
Assume that the unital completely positive map 
\[
\Phi_{\NN}\colon B\otimes C\ni\sum_{i=1}^na_i\otimes x_i
\mapsto \sum_{i=1}^n \ce(a_i)x_i\in\B(\hh)
\]
is continuous w.r.t.\ the minimal tensor norm on 
$B\otimes C$. 
Then $\NN$ is injective.
\end{lem}
\begin{proof}
Since $B\otimes_{\min}C\subset\B(\hh)\otimes_{\min}C$ 
and $\B(\hh)$ is injective, $\Phi_{\NN}$ extends to 
a unital completely positive map 
$\Psi\colon\B(\hh)\otimes_{\min}C\to\B(\hh)$. 
Then, $\Psi$ is automatically a $C$-bimodule map. 
Put $\psi(a)=\Psi(a\otimes1)$ for $a\in\B(\hh)$. 
Then, for every $a\in\B(\hh)$ and $x\in C$, we have 
\[
x\psi(a)=\Psi(1\otimes x)\Psi(a\otimes 1)
=\Psi(a\otimes x)=\Psi(a\otimes 1)\Psi(1\otimes x)=\psi(a)x.
\]
Hence, $\psi$ maps $\B(\hh)$ into $C'=\MM$. 
It follows that $\tilde{\psi}=\ce\psi\colon\B(\hh)\to\NN$ 
is a unital completely positive map such that $\tilde{\psi}|_B=\ce|_B$. 
Let $I$ be the set of all triples $(\EE,\FF,\e)$, 
where $\EE\subset\NN$ and $\FF\subset\NN_*$ are finite subsets 
and $\e>0$ is arbitrary. 
The set $I$ is directed by the order relation 
$(\EE_1,\FF_1,\e_1)\le(\EE_2,\FF_2,\e_2)$ 
if and only if $\EE_1\subset\EE_2$, $\FF_1\subset\FF_2$ and $\e_1\geq\e_2$.
Let $i=(\EE,\FF,\e)\in I$ and let $E\subset\NN$ be the finite dimensional 
operator system generated by $\EE$. 
We note that $E\subset\MM=pB^{**}$ and $\ce^*(\FF)\subset\MM_*=pB^*$, 
where $p\in B^{**}$ is the central projection supporting 
the identity representation of $B$ on $\hh$. 
Since $B$ is locally reflexive (cf.\ Lemma \ref{locref}), 
there is a unital completely positive map $\theta_i\colon E\to B$ 
such that for $a\in\EE$ and $f\in\FF$, we have 
\[
|\ip{\ce\theta_i(a),f}-\ip{a, f}|
=|\ip{\theta_i(a),\ce^*(f)}-\ip{a,\ce^*(f)}|<\e.
\]
Take a unital completely positive extension 
$\bar{\theta}_i\colon\B(\hh)\to\B(\hh)$ of $\theta_i$ 
and let $\sigma_i=\tilde{\psi}\bar{\theta}_i\colon\B(\hh)\to\NN$. 
It follows that for $a\in\EE$ and $f\in\FF$, we have 
\[
|\ip{\sigma_i(a),f}-\ip{a,f}|
=|\ip{\ce\theta_i(a),f}-\ip{a,f}|<\e.
\]
Therefore, any cluster point, in the point-ultraweak topology, 
of the net $\{\sigma_i\}_{i\in I}$ is a conditional expectation 
from $\B(\hh)$ onto $\NN$. 
\end{proof}
We now prove Theorem \ref{thm}, or more precisely, 
\begin{thm}
A finite von Neumann algebra $\MM$ satisfying 
the condition $\mathrm{(AO)}$ is solid.
\end{thm}
\begin{proof}
Let $\A$ be a diffuse von Neumann subalgebra in $\MM$. 
Passing to a subalgebra if necessary, we may assume $\A$ is abelian 
and prove the injectivity of $\NN=\A'\cap\MM$. 
It suffices to show $\Phi_{\NN}$ in Lemma \ref{mc} is continuous 
on $B\otimes_{\min}C$. 
Since $\A$ is diffuse, it is generated by a unitary $u\in\A$ 
such that $\lim_{k\to\infty}u^k=0$ ultraweakly.
Let $\Psi_n(a)=n^{-1}\sum_{k=1}^nu^kau^{-k}$ for $a\in\B(\hh)$ 
and let $\Psi\colon\B(\hh)\to\B(\hh)$ be its cluster point 
in the point-ultraweak topology. 
It is not hard to see that $\Psi|_{\MM}$ is a trace preserving 
conditional expectation onto $\NN$ and hence $\Psi|_{\MM}=\ce$. 
It follows that for any $\sum_{i=1}^na_i\otimes x_i\in\MM\otimes\MM'$, 
we have 
\[
\Psi(\sum_{i=1}^na_ix_i)=\sum_{i=1}^n\ce(a_i)x_i
=\Phi_{\NN}(\sum_{i=1}^na_i\otimes x_i).
\]
Since $\lim_{k\to\infty}u^k=0$ ultraweakly, 
we have $\K(\hh)\subset\ker\Psi$. 
This implies $\Psi=\tilde{\Psi}\pi$ for 
some unital completely positive map
$\tilde{\Psi}\colon\B(\hh)/\K(\hh)\to\B(\hh)$. 
Since $\nu$ in the condition~$\mathrm{(AO)}$ is continuous 
on $B\otimes_{\min}C$, so is $\Phi_{\NN}=\tilde{\Psi}\nu$.
\end{proof}
The following proposition was communicated to us by Popa \cite{popapc}. 
We are grateful to him for allowing us to present it here. 
\begin{prop}\label{popa}
Assume the type $\mathrm{II}_1$ factor $\MM$ (with separable predual)
contains a non-injective von Neumann subalgebra $\NN_0$ such that
$\NN_0'\cap\MM^\omega$ is a diffuse von Neumann algebra, where
$\MM^\omega$ is an ultrapower algebra of $\MM$.
Then there exists a non-injective von Neumann subalgebra
$\NN_1\subset\MM$ such that $\NN_1'\cap\MM$ is diffuse.
\end{prop} 
\begin{proof}
Replacing it with a subalgebra if necessary, 
we may assume the non-injective von Neumann subalgebra $\NN_0$ 
is generated by a finite set $\{x_1, x_2, \ldots, x_m\}$. 
By Connes' characterizations of injectivity 
(Theorem 5.1 in \cite{connes}),
it follows that there exists $\varepsilon>0$
such that if $\{x_1', \ldots, x_m'\}\subset\MM$ are
so that $\|x_i'-x_i\|_2 \leq \varepsilon$ then
$\{x_i'\}_i$ generates a non-injective von Neumann
subalgebra in $\MM$. 
Indeed, if there existed injective von Neumann algebras 
$\BB_k\subset\MM$ such that 
$\lim_k \| x_i - \e_{\BB_k}(x_i) \|_2=0,\ \forall i,$ 
then for any $\sum_{j=1}^n a_j\otimes y_j\in\NN_0\otimes\MM'$, 
we would have 
\begin{align*}
\| \sum_{j=1}^n a_j y_j\|_{\B(\hh)}
&\le\liminf_{k\to\infty}\| \sum_{j=1}^n \e_{\BB_k}(a_j) y_j\|_{\B(\hh)}\\
&=\liminf_{k\to\infty}\| \sum_{j=1}^n \e_{\BB_k}(a_j) \otimes y_j
 \|_{\BB_k\otimes_{\min}\MM'}
\le \| \sum_{j=1}^n a_j \otimes y_j\|_{\NN_0\otimes_{\min}\MM'},
\end{align*}
which would imply $C^*(\NN_0,\MM')\cong\NN_0\otimes_{\min}\MM'$ and 
the injectivity of $\NN_0$ (cf.\ the proof of Lemma \ref{mc}). 

Since $\NN_0 ' \cap \MM^{\omega}$ is diffuse, it follows
by induction that there exists a sequence
of mutually commuting, $\tau$-independent
two dimensional abelian $*$-subalgebras $\A_n\subset\MM$,
with minimal projections of trace $1/2$,
such that 
\[
\|\e_{\A_{n+1}'\cap\MM}(x_i) - x_i\|_2 < \varepsilon/2^{n+1},\ \forall i.
\]
But if we denote 
by $\BB_n = \A_1 \vee \A_2 \vee ... \vee \A_n$ 
then we also have
\[
\|\e_{\BB_n'\cap\MM}(x_i) - \e_{\BB_n'\cap\MM}(\e_{\A_{n+1}'\cap\MM}(x_i))\|_2
< \varepsilon/2^{n+1},\ \forall i.
\]
Since $\e_{\BB_n'\cap\MM}\circ \e_{\A_{n+1}'\cap\MM}
= \e_{\BB_{n+1}'\cap\MM}$, if we denote
by $\A=\vee_n \BB_n$ and take into account
that $\e_{\A'\cap\MM} = \lim_n \e_{\BB_n'\cap\MM}$
(see e.g., \cite{popak}), 
then by triangle inequalities we get
\[
\|x_i - \e_{\A'\cap\MM} (x_i)\|_2 \leq \varepsilon,\ \forall i.
\]
Thus, if we take $\NN_1$ to be the von Neumann algebra generated by
\[
x_i' = \e_{\A'\cap \MM} (x_i),\ 1\leq i \leq m,
\]
then  $\NN_1$ satisfies the required conditions.
\end{proof}

\section{Remark}
We remark the possibility that, 
for a hyperbolic group $\G$, the $*$-homomorphism 
$$\LG\otimes\crg\ni\sum_{i=1}^n a_i \otimes x_i 
 \mapsto \pi(\sum_{i=1}^na_i x_i) \in \B(\eg)/\K(\eg)$$ 
may be continuous w.r.t.\ the minimal tensor norm. 
If this is the case, then it would follow that 
a von Neumann subalgebra $\NN\subset\LG$ is injective 
if and only if 
\[
C^*(\NN,\crg)\cap\K(\eg)=\{0\},
\] 
which would reprove our results 
(modulo Theorem 2.1 in \cite{connes}). 
\medskip

\noindent\textbf{Acknowledgment.}
The author would like to thank Professor Sorin Popa for 
showing him Proposition \ref{popa}, and Professor Nigel Higson 
for providing him the references \cite{hg} and \cite{skandalis}.
This research was carried out while the author was visiting 
the University of California at Los Angeles under the support of 
the Japanese Society for the Promotion of Science
Postdoctoral Fellowships for Research Abroad.
He gratefully acknowledges the kind hospitality of UCLA. 

\end{document}